\documentclass[conference]{IEEEtran}
\IEEEoverridecommandlockouts

\usepackage[hidelinks]{hyperref}

\usepackage{cite}
\usepackage{amsmath,amssymb,amsfonts}
\usepackage{algorithmic}
\usepackage{graphicx}
\usepackage{textcomp}
\usepackage{xcolor}
\usepackage{siunitx}
\usepackage{bm}
\usepackage{tikz}
\usepackage{pgfplots}
\usepackage{comment}

\usetikzlibrary{calc}
\usetikzlibrary{fit}
\usetikzlibrary{shapes,arrows,positioning,calc,decorations.markings}
\tikzset{
  block/.style = {draw, fill=white, rectangle, minimum height=3em, minimum width=3em},
  tmp/.style  = {coordinate}, 
  sum/.style= {draw, fill=white, circle, node distance=1cm, inner sep = 0mm},
  branch/.style = {circle, fill=black,inner sep=0.35mm, outer sep = 0},
  input/.style = {coordinate},
  output/.style= {coordinate},
  pinstyle/.style = {pin edge={to-,thin,black}
  }
}
\def\BibTeX{{\rm B\kern-.05em{\sc i\kern-.025em b}\kern-.08em
    T\kern-.1667em\lower.7ex\hbox{E}\kern-.125emX}}
\pgfplotsset{compat=1.18}

\definecolor{y_color}{RGB}{60, 180, 75}
\definecolor{theta_color}{RGB}{245, 130, 48}
\definecolor{theta_dot_color}{RGB}{0, 130, 200}
\definecolor{y_dot_color}{RGB}{145, 30, 180}
\definecolor{bowling_ball_color}{RGB}{80, 207, 22}

\usetikzlibrary{angles,quotes} 

\begin{document}

\title{Closed-Loop Identification and Tracking Control of a Ballbot}

\author{\IEEEauthorblockN{1\textsuperscript{st} Tobias Fischer}
\IEEEauthorblockA{\textit{Institute for Engineering in Medicine} \\
\textit{Universität zu Lübeck}\\
Luebeck, Germany \\
https://orcid.org/0009-0006-2893-3277}
\and
\IEEEauthorblockN{2\textsuperscript{nd} Dimitrios S. Karachalios}
\IEEEauthorblockA{\textit{Institute for Engineering in Medicine} \\
\textit{Universität zu Lübeck}\\
Luebeck, Germany \\
dimitrios.karachalios@uni-luebeck.de}
\thanks{This work was funded by the Deutsche Forschungsgemeinschaft (DFG, German Research Foundation) - 419290163.}
\and
\IEEEauthorblockN{3\textsuperscript{rd} Ievgen Zhavzharov}
\IEEEauthorblockA{\textit{Institute for Engineering in Medicine} \\
\textit{Universität zu Lübeck}\\
Luebeck, Germany \\
ievgen.zhavzharov@uni-luebeck.de}
\and
\IEEEauthorblockN{4\textsuperscript{th} Hossam S. Abbas}
\IEEEauthorblockA{\textit{Institute for Engineering in Medicine} \\
\textit{Universität zu Lübeck}\\
Luebeck, Germany  \\
h.abbas@uni-luebeck.de}
}

\maketitle

\begin{abstract}
Identifying and controlling an unstable, underactuated robot to enable reference tracking is a challenging control problem. In this paper, a ballbot (robot balancing on a ball) is used as an experimental setup to demonstrate and test proposed strategies to tackle this control problem. A double-loop control system, including a state-feedback gain in the outer-loop and a Proportional-Integral-Derivative (PID) controller in the inner-loop, is presented to balance the system in its unstable equilibrium. Once stability is reached, the plant’s response to a designed excitation signal is measured and interpreted to identify the system’s dynamics. Hereby, the parameters of a linearized model of the ballbot are identified with prior knowledge about the structure of the nonlinear dynamics of the system. Based on an identified linear time-invariant (LTI) state-space model, a double-loop control strategy is considered to balance the real system and to allow reference tracking. A linear quadratic regulator (LQR) is designed offline and implemented in the inner-loop to ensure balance. In the outer-loop, the estimated dynamics forecast the system’s behavior online using a model-predictive-control (MPC) design to find the optimal control input for reference tracking. The experimental results demonstrate the applicability of the proposed strategies. 
\end{abstract}


\section{Introduction}\label{sec:introduction}
\subsection{Motivation}\label{subsec:motivation}
\begin{figure}[!t]
    \centering
    \includegraphics[scale=0.1]{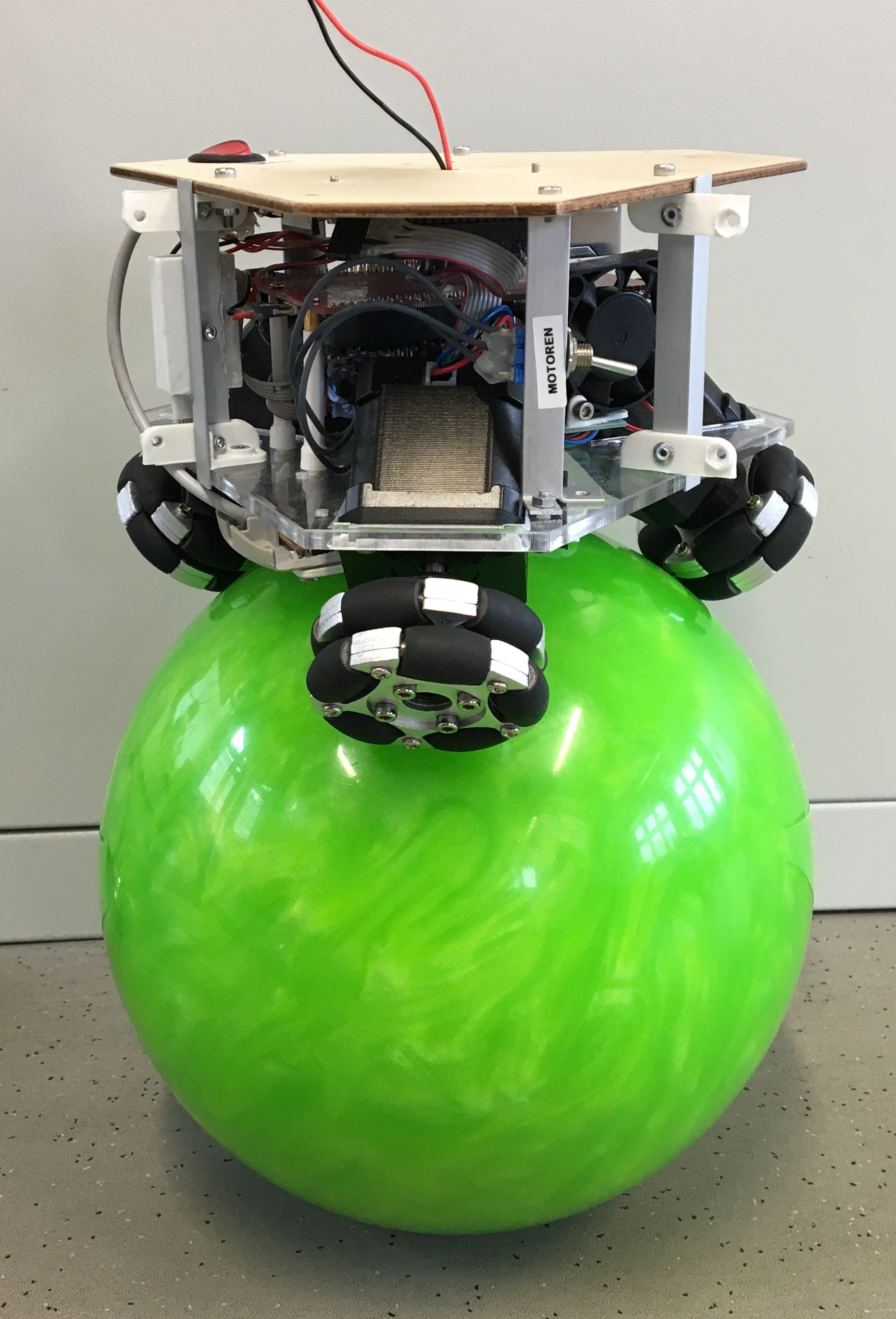}
    \caption{Picture of the ballbot (a side view) used in this work, which was built at the Institute for Electrical Engineering in Medicine, University of Lübeck.}
    \label{fig:picture-ime-ballbot-side}
\end{figure}

The ability of a ballbot, see Fig.~\ref{fig:picture-ime-ballbot-side}, to balance on its single spherical wheel distinguishes it from many other mobile robots since this grants omnidirectional movement. Achieving this capability requires the development of adequate control approaches. 
A ballbot provides five degrees of freedom in space. 
Since it performs movements by controlling a single ball, it is an under-actuated system that is challenging to control. However, a well-controlled ballbot can move in any desired direction in the horizontal plane without rotating around its vertical axis to align the body and path. This allows for a high degree of maneuverability in dynamically changing environments, such as the potential use case of service robotics. 

For most underactuated robotic systems, the dynamic model structure is well-known based on fundamental laws of physics. However, obtaining the physical parameters, particularly those related to inertia terms (components of the inertia tensor) and friction terms, requires specific information for each individual robot. These parameters can be obtained either through computations using classical mechanics \cite{SiScViOr08} or experimentally via direct measurements from the system. Note that the friction parameters must be determined experimentally. Identifying the physical parameters as lumped parameters usually provides real-world values related to the actual system. The computational approach can provide good model insight and is useful for concept studies and prototyping; however, it is time-consuming and demands a thorough understanding of the object. In this paper, we opt for the experimental approach.

The most challenging step in the experimental approach is stabilizing the system prior to identification to acquire the measurements for identification. Often, PID controllers are employed, and their constants are tuned using trial and error methods, see, e.g., \cite{PuBo20,Studt}, yielding the stabilization of the system being very cumbersome. To tackle this problem, we adopt the double-loop control approach proposed by \cite{Pham}, which offers meaningful insights into the stabilization problem and provides heuristics for tuning some of the controller gains.

The parameter identification based on closed-loop experiments using input/output measurements can result in highly biased parameters \cite{Van den Hof}, \cite{Ljung} when the input data are correlated with noisy output measurements. This is usually the case as the measured output is fed back to the controller, which calculates the input.
To address this issue, we adopt an indirect closed-loop identification approach \cite{Van den Hof}. In this approach, the identification is conducted using the output signal and an external signal that are uncorrelated with the noise contribution in the output signal. Consequently, unbiased estimates of the plant parameters can be obtained.

As a next step, model-based control strategies can be employed to control the system. For the ballbot, the ultimate goal is to enable the robot to navigate through cluttered environments and engage in human-robot interaction. Therefore, a promising candidate is model predictive control (MPC) \cite{MacIejowski}, which has the potential to address path-following tasks and trajectory tracking \ simultaneously cite{Jespersen} within its optimization problem. However, to fully utilize the capabilities of MPC, it is common practice, especially in controlling unmanned aerial vehicles \cite{Kamel}, to equip the system with an internal controller such as PID to maintain the stability of the closed-loop system while the MPC, operating in an outer loop, handles path-following or trajectory tracking tasks.

This paper presents a method for controlling an unknown, unstable, underactuated robot using the ballbot system as a testbed. To stabilize the system and acquire the necessary measurements from the real-time system, a double-loop approach inspired by \cite{Pham} is employed, which allows a multi-harmonic signal to excite the ballbot around its unstable equilibrium. Subsequently, an indirect closed-loop identification technique is utilized, leveraging the known structure of its linearized state-space dynamics to estimate its parameters. To support trajectory tracking on the real system, the identified state-space model is utilized in designing an LQR state feedback controller as an inner control loop to balance the ballbot. Furthermore, the identified model acts as a predictor in a model predictive control scheme configured to track a pre-specified reference signal. The reference tracking capability of the control approach is successfully validated on the physical system.



\subsection{Contents}\label{subsec:contents}
We start our analysis with some preliminaries in Sec.~\ref{sec:preliminaries}. 
consisting of the experimental setup 
and the mathematical model of the ballbot derived from first principles.
Further, 
Section~\ref{sec:parameter_identification} introduces the proposed identification framework 
including the double-loop approach 
to stabilize the system, 
perturbation of the plant with a multi-harmonic input that enables 
the required measurements for the identification task 
and the indirect closed-loop identification procedure 
resulting in the identified parameters. Finally, Sec.~\ref{sec:control_design} presents the control design procedure to
provide the control law that drives the physical ballbot to a given reference trajectory. The paper summarizes the results in Sec.~\ref{sec:control_design} and concludes in Sec.~\ref{sec:conclusion}.

\section{Preliminaries}\label{sec:preliminaries}

\subsection{Experimental Setup}\label{subsec:experimental_setup}
\begin{figure}[!t]
    \centering
    \includegraphics[width = 5cm, trim=500 50 200 100, clip, angle = 0]{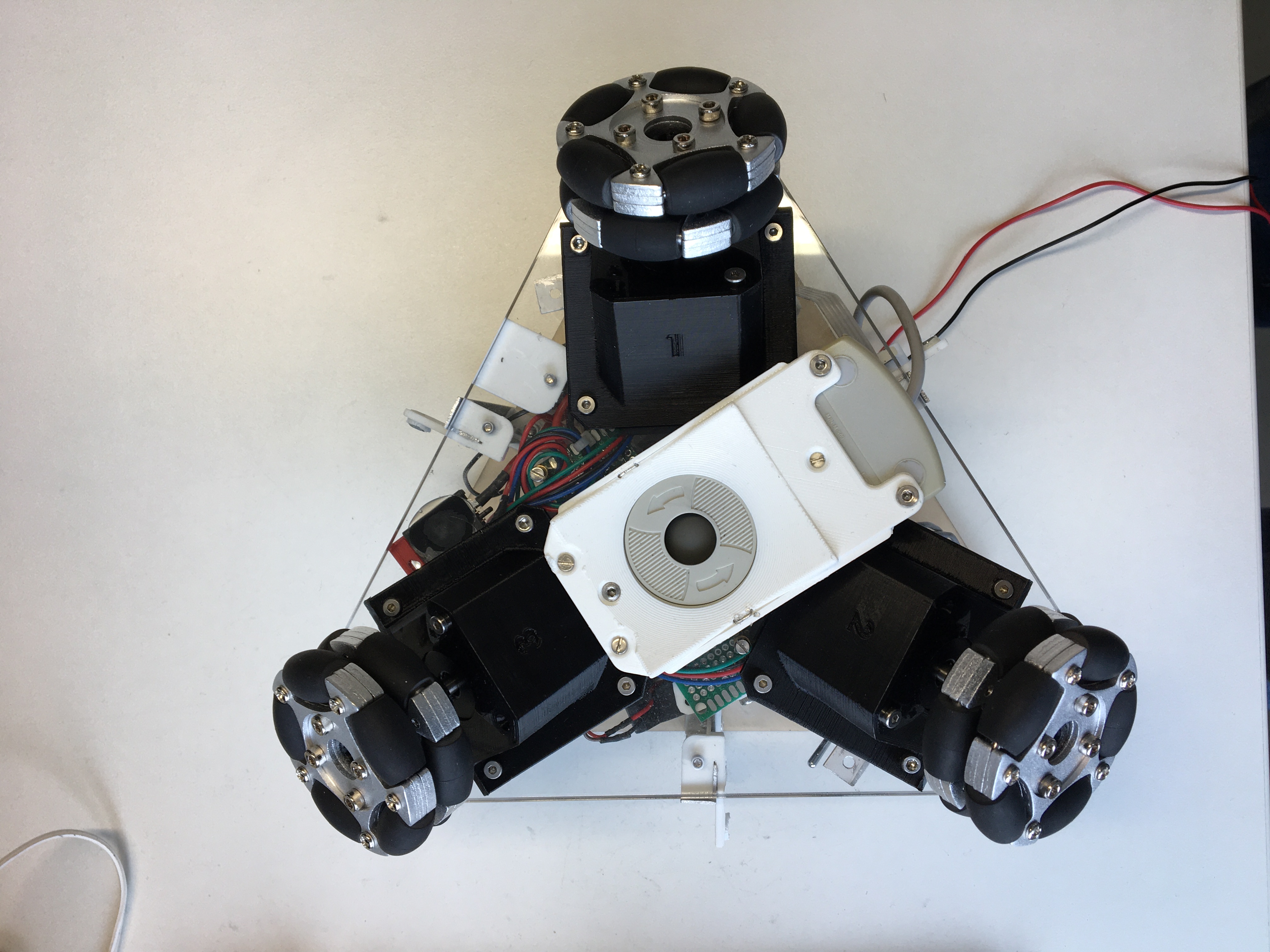}
    \caption{Picture of the IME ballbot, bottom view. The three motors and omniwheels form an equal triangle. The trackball sensor is positioned at its center, riding on top of the bowling ball.}
    \label{fig:picture-ime-ballbot-bottom}
\end{figure}

The ballbot is an experimental setup built at the Institute for Electrical Engineering in Medicine of the University of Lübeck. The foundation of the robot is a horizontal base plate to which the other components are attached. The robot’s actuators are three NEMA $17$ stepper motors operating at a voltage of $12$ (V) actuating three double-row omniwheels (with a diameter of 58 (mm)), which drive a bowling ball, on which the robot balances, see \autoref{fig:picture-ime-ballbot-side}.
An inertial measurement unit (IMU, SparkFun ICM30948) is used to provide information about the tilting angles of the ballbot. The inertia measurement unit (IMU) sampling time is set to $5$ (ms), as this will be the sampling time of the control algorithm. To determine the rotational speed of the bowling ball (relative to the robot itself), an adapted physical computer mouse (Mitsumi PS2) is mounted underneath the base plate of the ballbot, see \autoref{fig:picture-ime-ballbot-bottom}.

This mouse is referred to as the trackball sensor. It consists of the trackball itself running on the bowling ball’s surface and two smaller encoder wheels driven by the trackball. Two optical encoders are reading the rotational velocity of the encoder wheels towards the x and y directions, respectively.

A microcontroller (XMC4700 Relax Kit) is used for its capability to communicate with the hardware and its support for real-time operation, while the control algorithm runs on a small on-board computer (Raspberry Pi 4b 4GB) that provides a higher computing power. As a control input, the motors receive a velocity signal from the microcontroller regarding steps per second.

All experiments are conducted on a hard, low-friction surface. This means that even the ball alone is unstable, unlike on other surfaces like carpet.

\subsection{Mathematical Modeling}\label{sec:mathematical_modeling}
A coordinate system is placed at the center of the ball with the $z$-axis pointing upwards.
Two models of an identical structure \eqref{eq:nonlinear-model} are used to describe the behavior of the ballbot in space as done in \cite{Studt}. One model represents the behavior of the ballbot in the $xz$-plane, while the other reflects its behavior in the $yz$-plane. As an example, the dynamics of the $xz$-plane can be described  using first principles as follows:
\begin{subequations}
    \begin{align}
         M(q)\ddot q &+ \Tilde{C}(q, \dot q) + \Tilde{D}(\dot q) + G(q)=\Tilde{B}\tau_y,  \\
        q &=\begin{bmatrix}
            y & \theta_y
        \end{bmatrix}^T
    \end{align}
    \label{eq:nonlinear-model}
\end{subequations}
with the matrices 
\begin{align}
\begin{aligned}
\label{eq:parts}
    M &= \begin{bmatrix}
        b_{1} & -b_{2} + \ell r\cos{\theta_{y}} \\
        -b_{2} + \ell r\cos{\theta_{y}} & b_{3}
    \end{bmatrix}, \\
    \Tilde{C} &= \begin{bmatrix}
        -\ell r\sin({\theta_{y})\dot{\theta}_y^{2}} \\
        0
    \end{bmatrix},~\Tilde{D} = \begin{bmatrix}
        b_{4}\frac{\dot{y}}{r}\\[1mm]
        b_5\dot{\theta}_y
    \end{bmatrix}, \\
    G &= \begin{bmatrix}
        0 \\
        -\ell g\sin({\theta_{y}})
     \end{bmatrix},~\widetilde{B} = \begin{bmatrix}
        \frac{r}{r_{w}} \\
        - \frac{r}{r_{w}}
    \end{bmatrix},
\end{aligned}
\end{align}
where $\Tilde{B}, M, \Tilde{C}, \Tilde{D}, G$ represent the input matrix, the mass matrix, the vector of Coriolis and centripetal torques, the frictional torque vector, the vector of gravitational torques respectively, with  $b_1, b_2, \cdots, b_5, \ell$ represent combinations of dynamic and kinematic physical parameters of the system, $r$ is the radius of the ball,
$r_w$ is the radius of the omniwheels, $g$ is the gravitational
acceleration 
and $\tau_y$ is the torque that is applied by the motors in the $xz$-plane. 
The vector $q$ is the so-called generalized coordinates with $y$ indicating the position of the ball along the $y$-axis 
and $\theta_y$ 
the robot's tilt angle around the $y$-axis. 

The nonlinear model \eqref{eq:nonlinear-model} of the ballbot is linearized at the unstable equilibrium $q = 0,\dot q=0$ 
resulting in the continuous-time linear time-invariant (CT LTI) state-space representation  
\vspace{-4mm}
\begin{subequations}
    \begin{align}
        \dot x(t) &= Ax(t)+Bu(t)\\
        y(t)&=Cx(t)+Du(t),
    \end{align}
    \label{eq:ctlti system}
\end{subequations}%
\hspace{-2mm} where the particular model of the ballbot is given as
\begin{equation}
    \begin{bmatrix}
        \dot y\\
        \dot \theta_x\\
        \ddot y\\
        \ddot \theta_x
    \end{bmatrix} = \underbrace{\begin{bmatrix}
        0 & 0 & 1 & 0\\
        0 & 0 & 0 & 1\\
        0 & p_1 &  p_2 & p_7\\
        0 & p_4 & p_5 & p_8
    \end{bmatrix}}_{A}
    \begin{bmatrix}
        y\\
        \theta_x\\
        \dot y\\
        \dot \theta_x
    \end{bmatrix}+\underbrace{\begin{bmatrix}
        0\\ 0\\ p_3\\ p_6
    \end{bmatrix}}_{B}u_{\theta x}
    \label{eq:linearized-ballbot-model}
\end{equation}
 with the output matrix $C$ is considered as the identity matrix, as sensors can measure all system states,  the feed-through $D$ matrix is zero, and $p_1,p_2,\cdots,p_8$ are constant parameters which are related to the physical parameters in \eqref{eq:parts}.
The purpose of the closed-loop identification proposed in this work is to identify the values of these parameters.
 
\begin{figure}[!t]
    \newcommand{\myscale}{0.7}
    \centering
    \includegraphics[scale=0.8]{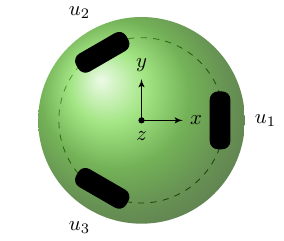}
    \caption{Schematic wheel arrangement, top view. The robot's wheels run on the bowling ball at the latitude angle of $\alpha = \si{45 \degree}$, denoted by the dashed line. As the motors do not apply torque directly in the $xz$ or the $yz$-planes, a power distribution \eqref{eq:planar_to_motors} is needed.}\label{fig:ballbot-wheel-arrangement}
\end{figure}

Since the controllers operate on the planar model of the ballbot, their control input also corresponds to the planar model. Therefore, it is important to note that the motors and omniwheels that drive the bowling ball are not in these planes, see \autoref{fig:ballbot-wheel-arrangement}. To distribute the control input of the planar model to the motors \eqref{eq:planar_to_motors} is used, where $\alpha = \SI{45}{\degree}$ is the angle of latitude where the omniwheels touch the ball. 
\begin{equation}
    \begin{bmatrix}
        u_1\\u_2\\u_3
    \end{bmatrix} = \frac{1}{3}
    \begin{bmatrix}
        \frac{2}{cos\alpha} & 0 & \frac{1}{sin\alpha}\\
        -\frac{1}{cos\alpha} & \frac{\sqrt{3}}{cos\alpha} & \frac{1}{sin\alpha}\\
        -\frac{1}{cos\alpha} & -\frac{\sqrt{3}}{cos\alpha} & \frac{1}{sin\alpha}        
    \end{bmatrix}
    \begin{bmatrix}
        u_{\theta x}\\u_{\theta y}\\u_{\theta z}
    \end{bmatrix}
    \label{eq:planar_to_motors}
\end{equation}

\section{Parameter Identification}\label{sec:parameter_identification}
As the given physical system is unstable, it is essential to stabilize it before identification. Thus, a balancing controller is designed and implemented. Next, a tuned perturbation signal is applied to this controlled system to collect a data set by which the model of the ballbot can be identified employing an 
excitation that can reveal the dynamics. The identified model will later be a part of the predictor within the MPC.
\subsection{Stabilizing the Unidentified System}\label{subsec:stabilizing_the_unidentified_system}
\begin{figure}[!t]
\centering
\includegraphics[scale=0.8]{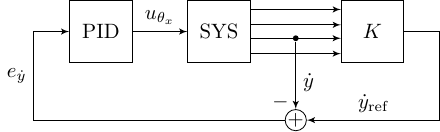}
\caption[Block diagram of stabilizing the unidentified ballbot model.]{Block diagram of the feedback controller used to stabilize the system. SYS represents the (unidentified yet measured) dynamics of the ballbot. Here, $K=[k_y\quad k_{\theta_x}\quad k_{\dot y}\quad k_{\dot\theta}]$.} \label{fig:pid-balancing-controller}
\end{figure}

All the states are measured by the sensors.
Digital filters are applied to the measured signals, and offset errors are removed. After the implementation of online signal processing, a controller to stabilize the system can be designed. However, this task is non-trivial since the system does exhibit nonminimum phase behavior and its dynamics are yet to be discovered. Therefore, a 
double-loop feedback control strategy inspired by \cite{Pham}  can balance the ballbot, see  \autoref{fig:pid-balancing-controller}, 
which mainly consists of weighted feedback of those three states relevant for balance. 


Assume that 
there  is a speed
$\dot y_{\text{ref}}$ 
of the bowling ball, which can stabilize its position; the task is to actuate the bowling ball to track this speed. The state feedback matrix $K$ represents the outer-loop controller, see  \autoref{fig:pid-balancing-controller}, that generates $\dot y_{\text{ref}}=Kx$ to keep the tilt angle $\theta$ at zero position. The difference $e_{\dot y} = \dot y_{\text{ref}} - \dot y$ between the desired speed and the actual speed is then fed back into a PID controller to actuate the ball accordingly, which represents the inner-loop controllers. Thus, the purpose of the PID controller is to allow $\dot y$ to follow $\dot y_{\text{ref}}$. The analysis of such a double-loop approach has been carried out in \cite{Pham}, where a  PI controller was considered in the inner loop. 

\begin{figure}[!t]
    \centering
    \includegraphics[scale=0.8]{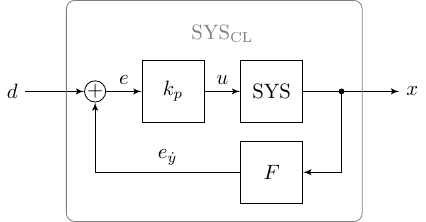}
    \caption{Block diagram of the control structure used for system identification. The excitation signal is the external input signal $d$. The feedback-matrix F denotes the same feedback structure as shown in \autoref{fig:pid-balancing-controller} to generate $e_{\dot y}$. This means that 
    $F = [0 \quad k_\theta \quad k_{\dot y}-1 \quad k_{\dot \theta}]$.}
    \label{fig:block-diagram-from-d-to-x}
\end{figure}

A discrete-time PID controller, \autoref{fig:pid-balancing-controller}, is implemented in real time as follows:
\begin{subequations}
    \begin{align}
        e_y(k) &= e_y(k-1) + e_{\dot y}(k) T_s,\\
        e_{\ddot y}(k) &= \frac{e_{\dot y}(k) - e_{\dot y}(k-1)}{T_s},\\
        u(k) &= K_P e_{\dot y}(k) + K_I e_y(k) + K_D e_{\ddot y}(k),
    \end{align}
    \label{eq:pid-implementation}
\end{subequations}
where 
$k$ denotes the discrete-time index, and $T_s$ is the sampling time of the control algorithm. The controller gains that stabilized the system are $K_P=180$, $K_I=830$, and $K_D=50$, respectively. 
Utilizing the PID controller enabled easy tuning of the feedback parameters $K=[k_y\quad k_{\theta_x}\quad k_{\dot y}\quad k_{\dot\theta}]$ to generate the reference speed.
\begin{table}[!t]
\caption[table head]{Parameter of the stabilizing controller in \autoref{fig:pid-balancing-controller} and \autoref{fig:block-diagram-mpc-lqr-control-structure}.}
    \centering
    \begin{tabular}{|l|c|c|c|c|c|}
        \hline
        Parameter & $k_p$ & $k_y$ & $k_{\theta_x}$ & $k_{\dot y}$ & $k_{\dot\theta_x}$\\
        \hline
        Value & 300 & 0 & 1.2 & 1.1 & 0.005\\
        \hline
    \end{tabular}
    \label{tab:parameters stabilizaion pid}
\end{table}
Once a robust balance is accomplished, simplifications of the control structure are considered to 
reduce the complexity of the subsequent 
indirect identification considered here, which requires the information of the controller. 
Therefore, replacing the PID with a P controller in the inner loop
still results in acceptable balancing performance. Hence, for the identification, the control structure is shown in \autoref{fig:block-diagram-from-d-to-x} is employed.
The parameters of all gains used for stabilization using the double-loop approach are presented in \autoref{tab:parameters stabilizaion pid}.
\subsection{The Perturbation Signal}\label{subsec:the_perturbation_signal}
The stabilized system can now be excited by a perturbation signal to 
acquire the necessary measurements 
to identify the state-space model of the ballbot. During the testing phase of the balancing controller, the system bandwidth of about \SI{1}{\hertz} could be observed, meaning the dominant movements were actuated at this rate/bandwidth. This delivers an indication of how to construct the perturbation signal to provide informative excitation of the system.
Utilizing this guess 
about the system’s behavior in the frequency domain, designing the excitation signal in frequency domain can also be considered an appropriate choice. Thus, the excitation signal is a composition of multiple sine waves at frequencies around the estimated bandwidth of the system. 
The signal used to perturb the ballbot, as shown in Fig.~\ref{fig:block-diagram-from-d-to-x}, is formulated as
\begin{equation}
    d(t) = \alpha\sum_{i = 1}^5 a_i \sin(b_i 2\pi t),
\label{eq:excitation-signal}
\end{equation}
where  the values of the parameters $ a_i, b_i$ are given  in \autoref{tab:parameters-excitation-signal}.
\begin{table}[!t]
\caption[table head]{Parameters of the excitation signal}
    \centering
    \begin{tabular}{ |c|c|c|c|c|c| }
        \hline
        $i$ & 1 & 2 & 3 & 4 & 5 \\
        \hline
        $a_i$ [\si{\centi\metre\per\second}] & 0.14 & 1 & 0.27 & 0.14 & 0.125\\
        \hline
        $b_i$ [\si{\hertz}]& 0.43 & 0.64 & 0.7 & 3.4 & 5.1\\
        \hline
    \end{tabular}
    \label{tab:parameters-excitation-signal}
\end{table}
In this way, the individual components of d are adjusted by changing $a_i$ and $b_i$, while
$\alpha$ is used to scale the overall amplitude of the excitation signal to not over-excite the ballbot. 
Increasing the value of $\alpha$ 
without destabilizing the system results in a wider range of excitation within the stable range. 
After collecting sufficient identification data, a model of the system can be estimated from these logged signals.

\subsection{Indirect Close-Loop Identification} 
\label{subsec:gray-box_identification}
Once the set of measurements including the external perturbation input $d$ and all the states of the ballbot as the output according to 
the closed-loop system in \autoref{fig:block-diagram-from-d-to-x} is collected, the system identification can be performed. 
This task aims to obtain the parameters of the linearized model \eqref{eq:linearized-ballbot-model}. We refer to this as an indirect closed-loop Identification as we identify the closed-loop system, while the information of the open-loop system is extracted given the controller's information. It is also possible to call it
gray box model identification, as a predefined model structure is enforced. It is done using the \textsc{Matlab} software based on the closed-loop system structure shown in \autoref{fig:block-diagram-from-d-to-x}. An important feature of this identification 
approach is that it is carried out using the external signal $d$, which is uncorrelated with the noise contribution in the feedback signals, i.e., the states, which guarantees an unbiased estimate of the parameters, \cite{Van den Hof}, \cite{Ljung}.


Since the state $y$ does not influence the balance of the ballbot ($k_y=0$),   it is not considered during the identification and will later be introduced to the model, as it can be easily derived from the state $\dot y$. The structure of the to-be-identified system SYS, see \autoref{fig:block-diagram-from-d-to-x}, is known, yet the collected data for system identification originates from the closed-loop system SYS$_\text{CL}$, \autoref{fig:block-diagram-from-d-to-x}.
Hence, the approach is to
\begin{enumerate}
 \item derive the \textit{structure} of the gray box model SYS$_{\text{CL}}$
 \item identify the parameters of the model SYS$_{\text{CL}}$,
 \item extract the parameters of SYS from SYS$_{\text{CL}}$.
\end{enumerate}
Combining SYS  and the feedback controller into one state-space system is achieved by formulating $u$ in terms of the state vector $x$ and the external input $d$ as shown in \eqref{eq:structure-of-sys-cl}. This yields SYS$_{\text{CL}}$, where the output matrix $C = C_{\text{CL}} = I_4$ and 
$D = D_{\text{CL}} = 0_{4 \times 1}$. 
Only the system matrix and the input matrix are altered, as they now also represent the feedback controller. Thus, the closed-loop system is represented by
\begin{subequations}
  \begin{align}
  \dot x &= \underbrace{\left (A + B k_pF\right )}_{A_{\text{CL}}} x + \underbrace{B k_p}_{B_{\text{CL}}}d\text{,}\\
 F&=\begin{bmatrix} 0 & k_\theta & k_{\dot y} - 1 & k_{\dot \theta}\end{bmatrix}.
    \end{align} \label{eq:structure-of-sys-cl}
  \end{subequations}

Since the structure of SYS$_{\text{CL}}$ is established, it can be implemented into the function \texttt{odefun}, and the parameters of this system can be estimated using the \textsc{Matlab} function \texttt{greyest(data,init\_sys)}. The next step is to verify the identification results.

\begin{figure}[!t]
    \centering
    \includegraphics[scale=1.2]{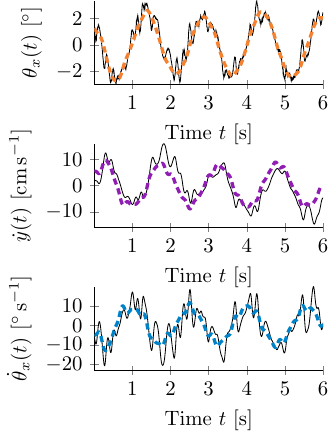}
    \caption{Output of the identified model SYS$_{\text{CL}}$ (dashed line) versus measured (solid line) signals based on the same external input signal $d$, see \autoref{fig:block-diagram-from-d-to-x}. The state $y$ is not shown here, as the position of the ballbot is not relevant for balancing.} 
    \label{fig:model-output}
\end{figure}

The quality of the estimated model are determined by
the fit rate (normalized root mean squared error)  calculated for the states $\theta_x$, $\dot y$, and $\dot \theta$ by 
\begin{equation}
    \text{fit} = 100\%\left (1-\frac{\|y - \hat y\|}{\|y - \text{mean}(y)\|}\right )\text{,}
\end{equation}
where $y$ is the measured system output, and $\hat y$ is the output of the simulated model.
The fit rates of the individual states are displayed in \autoref{tab:fit-rates}.
\begin{table}[!t]
\caption{Fit rates of the identified model SYS$_{\text{CL}}$}
    \centering
    \begin{tabular}{ |c|c|c|c|c|c|c|c|c| }
        \hline
        State & $\theta$ & $\dot y$ & $\dot \theta$\\
        \hline
        Fit [\si{\percent}] & 68.4885 & 45.6070 & 43.3884\\
        \hline
    \end{tabular}
    \label{tab:fit-rates}
\end{table}
A visual comparison of the physical system's behavior versus the behavior of the linearized simulated model, both being actuated by the same excitation signal, is presented in \autoref{fig:model-output}. Different sets of data are used for identification and validation. 

The identification result is an identified model of SYS$_{\text{CL}}$, from which the parameters of the actual model of the ballbot SYS can be extracted.
This is achieved by performing a comparison of the entries of $A_\text{CL}$ versus $A$, yielding the parameters $p_1$ to $p_8$ for \eqref{eq:linearized-ballbot-model}. Identification of the parameters $b_i,~i=1,\ldots,5$  in \eqref{eq:parts} that can make the nonlinear model accessible for nonlinear control can be done by utilizing nonlinear methods as in \cite{Karachalios}.

\begin{table}[!t]
\caption{Identfied parameters of the linearized model \eqref{eq:linearized-ballbot-model}, where $r$ is the radius of the bowling ball}
    \centering
    \begin{tabular}{ |@{\hspace{5pt}}c@{\hspace{5pt}}| *{8}{@{\hspace{5pt}}c@{\hspace{5pt}}|} }
    \hline
        $i$ & 1 & 2 & 3 & 4 & 5 & 6 & 7 & 8\\
        \hline
        $p_i$ & 0.25$r$ & -13.87$r$ & 0.01 & 3.95 & 214.68/$r$ & -0.28 & -6.05 & 5.50\\
        \hline
    \end{tabular}
    \label{tab:model-params}
\end{table}

Lastly, the state $y$ has to be introduced to the identified model, as this particular system state is used for reference tracking in the next section. Since $y$ can be obtained through integrating $\dot y$ (already a state of the system), adding an integrator to the model fulfills the task.
The resulting continuous-time state-space model has the structure planned in \eqref{eq:linearized-ballbot-model}. 
\section{Control Design and Experimental Results}\label{sec:control_design}
Before designing the tracking controller, 
we rebuild the balancing controller after the system has been identified. Subsequently, the predictor for the MPC can be constructed, validated, and implemented in real time. Finally, a reference signal can be applied. 

A Linear Quadratic Regulator (LQR) provides an optimal state feedback to stabilize the ballbot. To design and implement an LQR to operate in real time, the model should first be discretized, we consider  exact discretization using a zero-order hold and a sampling time of $T_s=0.005$ seconds, resulting in 
the discrete-time linear time-invariant (DT LTI) state space model:
\begin{subequations}
    \begin{align}
        x(k+1) &= A_dx(k)+B_du(k)\\
        y(k)&=C_dx(k)+D_du(k),
    \end{align}
    \label{eq:dt state space system}
\end{subequations}
where $A_d,B_d,C_d,D_d$ are the discrete-time system matrices.
Next, the feedback matrix $K_\text{LQR}$ is computed based on the discrete model and appropriate weighting matrices. This step is performed using the \textsc{Matlab} function \texttt{K\_LQR = dlqr(A, B, Q, R)}. The state weight is set to \texttt{Q = diag([20 100 10 50])} as this penalizes errors in $\theta_x$ (the second state) but allows for the necessary velocity $\dot y\neq 0$ (the third state) to stabilize the ballbot. The control weight is set to \texttt{R = 200}, whereby even larger changes in the control weight have no significant effect on the generated \texttt{K\_LQR}. This indicates that a specific amount of input energy has to be applied to stabilize the system regardless of its weight. 

\begin{figure}[!t]
    \centering
    \includegraphics[scale=0.8]{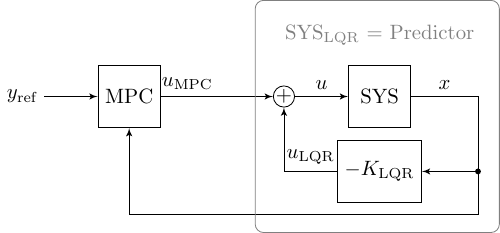}
    \caption{Block diagram of MPC LQR control structure}\label{fig:block-diagram-mpc-lqr-control-structure}
\end{figure}

The structure of the LQR balancing controller is shown in \autoref{fig:block-diagram-mpc-lqr-control-structure}. Feeding back the states through the generated feedback matrix $K_\text{LQR}$ stabilizes the system, as all eigenvalues of the discrete-time closed-loop system SYS$_\text{LQR}$ (see \autoref{fig:block-diagram-mpc-lqr-control-structure}) lie (strictly) within the unit circle.
Also, the LQR is tested in simulation, starting with nonzero initial conditions, safely controlling the system back to the equilibrium point.
\subsection{The Predictor}\label{sec:thepredictor}

The predictor applied in the MPC should not only represent the dynamics of the physical plant but should capture the behavior of the stabilized ballbot, including the LQR.
Consequently, the closed-loop system SYS$_\text{LQR}$ consisting of the identified model SYS of the ballbot as well as the feedback $-K_\text{LQR}x$, is used as a predictor in the MPC, see \autoref{fig:block-diagram-mpc-lqr-control-structure}.

Therefore, the system to be controlled is already internally stable (given no significant external disturbance), which prevents infeasibility that may arise from predicting the behavior of an unstable open loop system \cite{Cannon}.
This specific way of combining a predictive controller and an LQR is called the \textit{Dual-Mode Prediction Paradigm}.

The structure of the control approach for reference tracking is presented in \autoref{fig:block-diagram-mpc-lqr-control-structure}.
In addition to the predictor, implementing an MPC requires constraints, weights, and a prediction horizon. The MPC optimization problem applied here uses a quadratic cost function and some inequality constraints on the input and states of the ballbot:
\begin{subequations}
    \begin{align}
        \min_{u_0,\cdots,u_{N-1}}& e_N^TQ_Ne_N+\sum_{k=0}^{N-1}e_k^TQe_k+u_k^TRu_k\\
        \text{subject to }
        & e_k=x_k-x^{\rm r}_k\\
        &x_{k+1} = A_dx_k+B_du_k\\
        &|\theta_x| \leq 3\si{\degree}\\
        &|\dot y| \leq 15\text{ \si{\centi\metre\per\second}}\\
        &|\dot \theta_x| \leq 25\text{\si{\degree\per\second}}\\
        &|u_{\text{MPC}}| \leq 1000,
    \end{align}
\end{subequations}
where, $x^{\rm r}_k$ is the reference value for the ballbot states at discrete time instant $k$.
The reference tracking is be performed in the $y$-direction so that a path in space can later be handed to the MPC to follow. The weights are tuned and tested in simulation to 
\begin{subequations}
\begin{align}
    Q_N = Q &= diag(10^3,0,0,0)\text{,}\\
    R &= 0.3\text{.}
\end{align}
\end{subequations}
The ballbot's dynamics are to be considered to determine the duration of the prediction horizon.
Balancing the robot in its initial position, the control actions do not require online prediction since they depend linearly on the system's states.
However, regarding reference tracking, the system's special dynamic properties require a certain prediction.
Assuming a steady state, to achieve a particular velocity $\dot y$, the ballbot has to reach and maintain a corresponding tilt angle $\theta_x$, automatically resulting in a velocity of $\dot y \neq 0$ \si{\centi\metre\per\second}. Reaching the tilt angle requires actuation opposite to the actuation needed to maintain the reached angle.
Similar behavior is typical for an inverted pendulum, being a non-minimum-phase system.
The duration of these phases should be exceeded by the prediction horizon of the MPC.
Consequently, simulations are performed to determine this period.
The test scenario of a $20$ \si{\centi\metre} step in $y_{\text{ref}}$ shows that a prediction window of 4 \si{\second} is sufficient to enable reference tracking in simulation.
\subsection{Real-Time Implementation} \label{sec:mpc_real-time-implementation}
For an MPC to run on a small single-board computer (Raspberry PI) online, its computational demand should be minimized. This can be achieved by reducing the prediction window or increasing the simulation sample-time. Since the prediction window can not be arbitrarily short, the sample-time of the MPC is increased to $0.1$ \si{\second} while the LQR still operates with a sample-time of $0.005$ \si{\second}.

The library \cite{osqp-eigen}, which implements the solver \cite{osqp}, comes with an implementation of an MPC algorithm for simulation. This software is customized to work in real time and to suit the chosen predictor. Further, it is adapted to accept a dynamic reference value instead of a constant.

\begin{figure}[!t]
    \centering
  \includegraphics[scale=1.2]{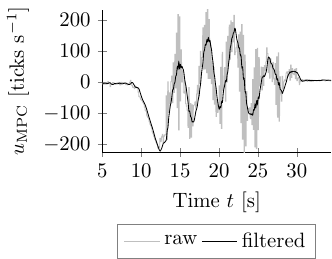}
    \caption{The raw and filtered version of $u_{\text{MPC}}$. The plotted information is from the same experiment as in \autoref{fig:tracking-result-smooth-step}. The unit of $u_{\text{MPC}}$ is ticks per second, as this signal is then together with $u_{\text{LQR}}$ using \eqref{eq:planar_to_motors} passed to the stepper motors (with micro stepping factor 4, they have 800 ticks per revolution).}
    \label{fig:mpc-output-raw-vs-filterd}
\end{figure}

Due to the MPC's susceptibility to noise and to provide a smooth control input $u_{\text{MPC}}$ to the system (see \autoref{fig:block-diagram-mpc-lqr-control-structure}), a second-order Butterworth low-pass filter is applied to the output of the MPC. The cutoff frequency is set to $1$ \si{\hertz}, and the raw and filtered versions of $u_{\text{MPC}}$ are plotted in \autoref{fig:mpc-output-raw-vs-filterd}. The noise of a higher amplitude appears when the ballbot dynamically tracks a reference.
Although the simulation of the MPC indicates the configuration of the parameters of the MPC, they need to be fine-tuned to match the behavior of the nonlinear physical system. Adjustment requires iterative testing while manually altering the weights, constraints, and sampling time.
\subsection{Reference Tracking Results}\label{sec:reference}
\begin{figure}[!t]
    \centering
    \includegraphics[scale=1.2]{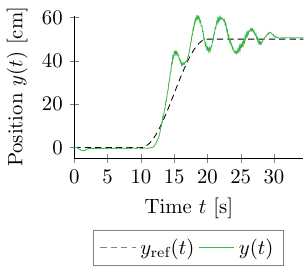}
    \caption{Plot of the experimental tracking result in $y$-direction over time. In the transitioning phase, the system leaves its linearization point and loses considerable robustness. However, it eventually regains a calm balance.}
    \label{fig:tracking-result-smooth-step}
\end{figure}

Once the MPC is set up to work in real time, the next step is to provide a reference trajectory. A smooth step signal is generated using a step signal with a sinusoidal transition.
This signal and the corresponding experimental result are plotted in \autoref{fig:tracking-result-smooth-step}.

\begin{figure}[!t]
    \centering
    \includegraphics[scale=1.2]{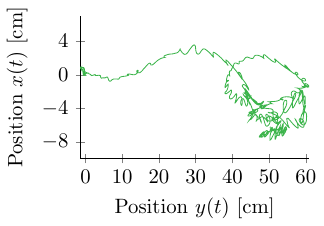}
    \caption{Plot of the experimental tracking result as a path in the $xy$-plane. Here, the $y$-position of the ballbot is controlled using the double-loop approach from \autoref{fig:block-diagram-mpc-lqr-control-structure} and the $x$-position is not controlled at all. The plotted information is from the same experiment as in \autoref{fig:tracking-result-smooth-step}.}
    \label{fig:tracking-result-smooth-step-xy-graph}
\end{figure}



Now the LQR introduced above is combined 
with the MPC to enable the ballbot to follow the provided reference trajectory.
In this double loop approach,  the LQR remains responsible for balancing the ballbot while the MPC adds station keeping capability (see \autoref{fig:block-diagram-mpc-lqr-control-structure} for the control structure).
As a result, the ballbot is stable in all states of the model.

Reference tracking, as the final objective of this work, is performed in the $yz$-plane with a time-dependent reference signal (trajectory) determining the position of the ballbot on the $y$-axis.
Reliable reference tracking is achieved in simulation, followed by experimental reference tracking utilizing the physical plant.
The physical ballbot is successfully controlled to transition to a new position where it regains a stationary balance (see \autoref{fig:tracking-result-smooth-step}).
However, during the transition phase, the experimental setup experiences significant disturbances, particularly during the initial jump of the trajectory. This indicates a high gain of the closed-loop system within this range of disturbance frequencies, resulting in the amplification of their amplitudes. 
Nevertheless, the reference tracking was successful, \autoref{fig:tracking-result-smooth-step}, and the ballbot settles after about 10 seconds.  This work demonstrates the feasibility of reference tracking using a linear predictor to control the experimental ballbot.

\section{Conclusion}\label{sec:conclusion}

This paper elaborates on a comprehensive concept from a given physical ballbot to its first reference tracking results. The methods can be applied to similar control problems, especially those involving plants requiring control at an unstable equilibrium. The stabilization of the ballbot has been achieved in three ways. Firstly, hand-tuned state feedback in combination with a PID controller is implemented, resulting in robust balance, thus enabling the indirect closed-loop system identification to obtain a linear model of the ballbot, which is needed for the model-based control approaches. The high quality of the identified linear state-space model determined from the real-time measurements demonstrates the approach's applicability.

Secondly, optimal state feedback utilizing an LQR based on the identified model has been successfully applied to the ballbot. As the first two approaches do not use the displacement distance for feedback, the ballbot is balanced but not controlled to stay at its initial position. Therefore, the third control approach combines the LQR used in the second approach with the  MPC. The LQR balances the ballbot while the MPC adds station-keeping capability in this double-loop approach. As a result, the ballbot is stable in all states of the model and can achieve successful tracking of a smoothed step reference trajectory. Continuing previous related studies, this work introduces the next step toward an autonomous ballbot, where path planning will be integrated into future endeavors.


\section*{Acknowledgment}
The work has been carried out at the Institute for Electrical Engineering in Medicine at the University of Luebeck. Conflict of interest: Authors state no conflict of interest.

\end{document}